\pdfoutput=1
\documentclass{amsart}
\usepackage[utf8]{inputenc}

\usepackage[margin=1in]{geometry}
\usepackage[expansion=false]{microtype}

\usepackage{amsmath, amsfonts, amssymb, amsthm, amsopn}
\usepackage{eucal}

\usepackage{tikz}
\usetikzlibrary{diagrams, decorations.pathreplacing}

\usepackage[pdfusetitle,unicode,hidelinks,bookmarksopen]{hyperref}

\newcommand{\resp}{{\sfcode`\.1000 resp.}}
\newcommand{\ie}{{\sfcode`\.1000 i.e.}}
\newcommand{\eg}{{\sfcode`\.1000 e.g.}}

\numberwithin{equation}{section}

\theoremstyle{plain}
\newtheorem*{theorem*}{Theorem}

\newtheorem{theorem}[equation]{Theorem}

\newtheorem{proposition}[equation]{Proposition}
\newtheorem{lemma}[equation]{Lemma}
\newtheorem{corollary}[equation]{Corollary}

\theoremstyle{definition}
\newtheorem{definition}[equation]{Definition}

\newtheorem{example}[equation]{Example}

\newtheorem{remark}[equation]{Remark}

\let\scr=\mathcal
\let\frak=\mathfrak
\let\bb=\mathbb
\def\N{\bb N}
\def\Z{\bb Z}
\def\Q{\bb Q}

\def\A{\bb A}
\def\P{\bb P}
\def\V{\bb V}
\def\1{\mathbf 1}

\def\G{\mathbb G}
\def\ph{\mathord-}
\def\pt{{\scriptscriptstyle\bullet}}

\let\del=\partial

\let\into=\hookrightarrow
\let\onto=\twoheadrightarrow
\let\tens=\otimes

\def\on{\textnormal{ on }}

\DeclareMathOperator{\Sym}{Sym}
\def\id{\mathrm{id}}

\DeclareMathOperator{\Hom}{Hom}
\DeclareMathOperator{\End}{End}

\DeclareMathOperator{\Map}{Map{}}
\DeclareMathOperator{\Mor}{Map^{Sp}{}}

\DeclareMathOperator{\Spec}{Spec}
\DeclareMathOperator{\Proj}{Proj}

\def\Nis{\mathrm{Nis}}

\def\Zar{\mathrm{Zar}}
\def\mot{\mathrm{mot}}

\def\fppf{\mathrm{fppf}}

\def\abs#1{\lvert #1\rvert}

\def\QCoh{\mathrm{QCoh}{}}

\let\cat=\mathrm

\def\Gr{\mathrm{Gr}{}}

\def\Sph{\cat S\mathrm{ph}{}}
\def\Sp{\cat S\mathrm{p}}

\def\H{\cat H}
\def\SH{\cat S\cat H}
\def\Mod{\cat{M}\mathrm{od}{}}
\def\Sch{\cat{S}\mathrm{ch}{}}

\def\Sm{{\cat{S}\mathrm{m}}}

\def\op{\mathrm{op}}

\def\Vect{\mathrm{Vect}{}}
\def\htp{\mathrm{htp}}

\def\KH{K\mskip -1mu H}

\def\Perf{\mathrm{Perf}{}}

\def\Pr{\mathcal{P}\mathrm{r}}
\def\CAlg{\mathrm{CAlg}}
\def\CMon{\mathrm{CMon}}

\def\minus{\smallsetminus}

\def\Cat{\mathcal{C}\mathrm{at}{}}
\def\PSh{\mathcal{P}}

\def\GL{\mathrm{GL}}

\def\KGL{\mathrm{KGL}{}}
\def\Fun{\mathrm{Fun}}
\def\fin{\mathrm{fin}}

\def\Tame{\mathrm{tqStk}{}}

\DeclareMathOperator{\Stab}{Stab}
\def\lax{\mathrm{lax}}

\def\sh{\mathrm{sh}}

\def\BG{{\mathbf BG}}

\let\lim=\relax
\DeclareMathOperator*{\lim}{lim}

\DeclareMathOperator*{\colim}{colim}

\title{Cdh descent in equivariant homotopy $K$-theory}
\author{Marc Hoyois}
\date{\today}
\thanks{The author was partially supported by NSF grant DMS-1508096}
\address{Fakultät für Mathematik, Universität Regensburg, Universitätsstr. 31, 93040 Regensburg, Germany}
\email{\href{mailto:marc.hoyois@ur.de}{marc.hoyois@ur.de}}
\urladdr{\url{http://www.mathematik.ur.de/hoyois/}}
\begin{document}

\begin{abstract}
We construct geometric models for classifying spaces of linear algebraic groups in $G$-equivariant motivic homotopy theory, where $G$ is a tame group scheme. As a consequence, we show that the equivariant motivic spectrum representing the homotopy $K$-theory of $G$-schemes (which we construct as an $E_\infty$-ring) is stable under arbitrary base change, and we deduce that the homotopy $K$-theory of $G$-schemes satisfies cdh descent.
\end{abstract}

\maketitle
\tableofcontents
\section{Introduction}

Let $K(X)$ and $K^B(X)$ denote the connective and nonconnective $K$-theory spectra of a quasi-compact quasi-separated scheme $X$ \cite{TT}.
The homotopy $K$-theory spectrum $\KH(X)$ was introduced by Weibel in \cite{Weibel}: it is the geometric realization of the simplicial spectrum $K^B(\Delta^\bullet\times X)$, where
\[\textstyle\Delta^n=\Spec\Z[t_0,\dotsc,t_n]/(\sum_i t_i-1)\]
is the standard algebraic $n$-simplex. There are natural transformations $K\to K^B\to \KH$, which are equivalences on regular schemes. 

Haesemeyer \cite{Haesemeyer} (in characteristic zero) and Cisinski \cite{Cisinski} (in general) proved that homotopy $K$-theory satisfies descent for Voevodsky's cdh topology. This was a key ingredient in the proof of Weibel's vanishing conjecture for negative $K$-theory, established in characteristic zero by Cortiñas, Haesemeyer, Schlichting, and Weibel \cite{CHSW}, and up to $p$-torsion in characteristic $p>0$ by Kelly \cite{Kelly} (with a simplified proof by Kerz and Strunk \cite{KerzStrunk}). More recently, Kerz, Strunk, and Tamme proved that $K$-theory satisfies ``pro-cdh descent'' and deduced Weibel's conjecture in complete generality \cite{KST}.

The goal of this paper is to extend the cdh descent result of Cisinski to a suitable class of Artin stacks, namely, quotients of schemes by linearizable actions of linearly reductive algebraic groups.
We will introduce a reasonable definition of the homotopy $K$-theory spectrum $\KH(\frak X)$ for such a stack $\frak X$, which agrees with $K(\frak X)$ when $\frak X$ is regular. The ``obvious'' extension of Weibel's definition works well for quotients by finite or diagonalizable groups, but, for reasons we will explain below, a more complicated definition is preferred in general.
Our main results are summarized in Theorem~\ref{thm:intro} below. In a sequel to this paper, joint with Amalendu Krishna, we use these results to prove vanishing theorems for the negative $K$-theory of tame Artin stacks \cite{HK}.

Let us first introduce some terminology.
A morphism of stacks $\frak Y\to\frak X$ will be called \emph{quasi-projective} if there exists a finitely generated quasi-coherent module $\scr E$ over $\frak X$ and a quasi-compact immersion $\frak Y\into \P(\scr E)$ over $\frak X$.
We say that a stack $\frak X$ has the \emph{resolution property} if every finitely generated quasi-coherent module over $\frak X$ is the quotient of a locally free module of finite rank.
Throughout this paper, we will work over a fixed quasi-compact separated (qcs) base scheme $B$, and we will say that a morphism of $B$-stacks $\frak Y\to\frak X$ is \emph{N-quasi-projective} if it is quasi-projective Nisnevich-locally on $B$.
We refer to \cite[\S2.7]{Hoyois} for the precise definition of a \emph{tame} group scheme over $B$. The main examples of interest are:
\begin{itemize}
	\item finite locally free groups of order invertible on $B$;
	\item groups of multiplicative type;
	\item reductive groups, if $B$ has characteristic $0$ (\ie, there exists $B\to \Spec\Q$).
\end{itemize}
Let $\Tame_B$ denote the $2$-category of finitely presented $B$-stacks that have the resolution property, that are global quotient stacks $[X/G]$ for some tame affine group scheme $G$, and such that the resulting map $[X/G]\to\mathbf BG$ is N-quasi-projective.\footnote{Quasi-projective $\mathbf BG$-stacks almost have the resolution property: they admit a schematic Nisnevich cover whose Čech nerve consists entirely of stacks with the resolution property \cite[Lemma 3.11]{Hoyois}. In particular, the requirement that stacks in $\Tame_B$ have the resolution property is immaterial as far as Nisnevich sheaves are concerned. Allowing $\mathbf BG$-stacks that are only \emph{Nisnevich-locally} quasi-projective is necessary to dispense with isotriviality conditions on $G$ when the base $B$ is not geometrically unibranch (see \cite[Remark 2.9]{Hoyois}).}
 For $\frak X\in\Tame_B$, we let $\Sch_{\frak X}\subset(\Tame_B)_{/\frak X}$ be the full subcategory of N-quasi-projective $\frak X$-stacks. 
The Nisnevich (\resp{} cdh) topology on $\Sch_\frak X$ is as usual the Grothendieck topology generated by Nisnevich squares (\resp{} Nisnevich squares and abstract blowup squares). The Nisnevich and cdh topologies on $\Tame_B$ are generated by the corresponding topologies on the slices $\Sch_\frak X$.

\begin{remark}
	If $B$ has characteristic zero, the $2$-category $\Tame_B$ includes all Artin stacks of finite presentation, with affine stabilizers, and satisfying the resolution property. Indeed, by a theorem of Gross \cite[Theorem 1.1]{Gross}, such stacks have the form $[X/\GL_n]$, where $X$ is a quasi-affine $\GL_n$-scheme. 
\end{remark}

\begin{remark}
	The stacks in $\Tame_B$ share many features with the ``tame Artin stacks'' considered in \cite{AOV}. There are two essential differences: our stacks are not required to have finite diagonal, but theirs are not required to have the resolution property. 
\end{remark}

\begin{theorem}\label{thm:intro}
	Let $B$ be a quasi-compact separated base scheme.
	 There exists a cdh sheaf of $E_\infty$-ring spectra 
	 \[\KH\colon \Tame_B^\op\to\CAlg(\Sp)\]
	 and an $E_\infty$-map $K\to \KH$ with the following properties.
	\begin{enumerate}
		\item If $\frak X\in\Tame_B$ is regular, the map $K(\frak X)\to \KH(\frak X)$ is an equivalence.
		\item $\KH$ is homotopy invariant in the following strong sense: if $p\colon\frak Y\to\frak X$ is an fpqc torsor under a vector bundle, then $p^*\colon \KH(\frak X)\to\KH(\frak Y)$ is an equivalence.
		\item $\KH$ satisfies Bott periodicity: for every vector bundle $\frak V$ over $\frak X$, there is a canonical equivalence of $\KH(\frak X)$-modules $\KH(\frak V\on\frak X)\simeq\KH(\frak X)$.
		\item Suppose that $\frak X\in\Sch_{\mathbf BG}$ where $G$ is an extension of a finite group scheme by a Nisnevich-locally diagonalizable group scheme.
		 Then $\KH(\frak X)$ is the geometric realization of the simplicial spectrum $K^B(\Delta^\bullet\times\frak X)$.
	\end{enumerate}
\end{theorem}

From property (1) and the hypercompleteness of the cdh topology, we immediately deduce:

\begin{corollary}
	Suppose that $B$ is noetherian of finite Krull dimension and that every stack in $\Tame_B$ admits a cdh cover by regular stacks, \eg, $B$ is essentially of finite type over a field of characteristic zero. Then the canonical map $K\to \KH$ exhibits $\KH$ as the cdh sheafification of $K$.
\end{corollary}

The fact that $\KH$ is a cdh sheaf means that it is a Nisnevich sheaf and that, for every cartesian square
	\begin{tikzmath}
		\diagram{
		\frak W & \frak Y \\ \frak Z & \frak X \\
		};
		\arrows (11-) edge[c->] (-12)  (11) edge (21) (21-) edge[c->] node[below]{$i$} (-22) (12) edge node[right]{$p$} (22);
	\end{tikzmath}
in $\Tame_B$ such that $i$ is a closed immersion, $p$ is N-projective, and $p$ induces an isomorphism $\frak Y\minus\frak W\simeq\frak X\minus\frak Z$, the induced square of spectra
	\begin{tikzmath}
		\def\colsep{2em}
		\diagram{
		\KH(\frak X) & \KH(\frak Z) \\
		\KH(\frak Y) & \KH(\frak W) \\
		};
		\arrows (11-) edge (-12)  (11) edge (21) (21-) edge (-22) (12) edge (22);
	\end{tikzmath}
is cartesian. For quotients by finite discrete groups, we can improve this result as follows:

\begin{theorem}\label{thm:intro-2}
	Let $G$ be a finite discrete group and let
	\begin{tikzmath}
		\diagram{
		W & Y \\ Z & X \\
		};
		\arrows (11-) edge[c->] (-12)  (11) edge (21) (21-) edge[c->] node[below]{$i$} (-22) (12) edge node[right]{$p$} (22);
	\end{tikzmath}
	be a cartesian square of locally affine qcs $G$-schemes over $\Z[1/\abs G]$, where $i$ is a closed immersion, $p$ is proper, and $p$ induces an isomorphism $Y\minus W\simeq X\minus Z$. Then the induced square of spectra
	\begin{tikzmath}
		\def\colsep{2em}
		\diagram{
		\KH([X/G]) & \KH([Z/G]) \\
		\KH([Y/G]) & \KH([W/G]) \\
		};
		\arrows (11-) edge (-12)  (11) edge (21) (21-) edge (-22) (12) edge (22);
	\end{tikzmath}
	is cartesian, where $\KH(\frak X)$ denotes the geometric realization of the simplicial spectrum $K^B(\Delta^\bullet\times\frak X)$.
\end{theorem}

\begin{remark}
	If $G$ is a finite discrete group acting on a qcs scheme $X$, then $X$ is a locally affine $G$-scheme if and only if the coarse moduli space of the Deligne–Mumford stack $[X/G]$ is a scheme \cite[Remark 4.5]{RydhCMS}.
\end{remark}

We make a few comments on homotopy invariance.
As we observed in \cite{Hoyois},
most of the interesting properties of homotopy invariant Nisnevich sheaves on schemes only extend to stacks if homotopy invariance is understood in the strong sense of property (2) of Theorem~\ref{thm:intro}.
A typical example of a homotopy equivalence in that sense
is the quotient map $\frak X\to\frak X/U$ where $U$ is a split unipotent group acting on $\frak X$; this map is usually not an $\A^1$-homotopy equivalence, not even Nisnevich-locally on the target. This explains why our definition of $\KH$ for general stacks is more complicated than it is for schemes. Property (4) of Theorem~\ref{thm:intro} is explained by the fact that vector bundle torsors over such stacks are Nisnevich-locally split.

Properties (1)--(4) of Theorem~\ref{thm:intro} will essentially be enforced by the definition of the homotopy $K$-theory presheaf $\KH$ and foundational results on equivariant $K$-theory due to Thomason \cite{ThomasonK} and Krishna–Ravi \cite{KrishnaRavi}. The content of Theorem~\ref{thm:intro} is thus the statement that $\KH$ is a cdh sheaf.
Its proof uses the machinery of stable equivariant motivic homotopy theory developed in \cite{Hoyois}. Namely, the fact that $\KH$ is a Nisnevich sheaf satisfying properties (2) and (3) of Theorem~\ref{thm:intro} implies that its restriction to smooth N-quasi-projective $\frak X$-stacks is representable by a motivic spectrum $\KGL_{\frak X}\in\SH(\frak X)$. By \cite[Corollary 6.25]{Hoyois}, we can then deduce that $\KH$ satisfies cdh descent, provided that the family of motivic spectra $\{\KGL_{\frak X}\}_{\frak X\in\Tame_B}$ is stable under N-quasi-projective base change. This base change property is thus the heart of the proof. We will verify it by adapting Morel and Voevodsky's geometric construction of classifying spaces \cite[\S4.2]{MV} to the equivariant setting.

\begin{theorem}\label{thm:intro-3}
	For every $\frak X\in\Tame_B$, there exists an $E_\infty$-algebra $\KGL_{\frak X}\in\SH(\frak X)$ representing the $E_\infty$-ring-valued presheaf $\KH$ on smooth N-quasi-projective $\frak X$-stacks.
	Moreover, the assignment $\frak X\mapsto \KGL_{\frak X}$ is a section of $\CAlg(\SH(\ph))$ over $\Tame_B^\op$ that is cocartesian over N-quasi-projective morphisms. In particular, for $f\colon\frak Y\to\frak X$ N-quasi-projective, $f^*(\KGL_{\frak X})\simeq \KGL_{\frak Y}$.
\end{theorem}

Finally, we will observe that the Borel–Moore homology theory represented by $\KGL_{\frak X}$ on N-quasi-projective $\frak X$-stacks, for $\frak X$ regular, is the $K$-theory of coherent sheaves, also known as $G$-theory.

\begin{remark}\label{rmk:KR}
	In the paper \cite{KrishnaRavi}, the authors work over a base field. This assumption is used via \cite{HallRydh2} to ensure that the $\infty$-category $\QCoh(\frak X)$ of quasi-coherent sheaves is compactly generated and that the structure sheaf $\scr O_\frak X$ is compact.
	We claim that this holds for any $\frak X\in\Tame_B$.
	If $G$ is a linearly reductive affine group scheme, then $\scr O_{\mathbf BG}$ is compact in $\QCoh(\mathbf BG)$, by \cite[Theorem C (3)$\Rightarrow$(1)]{HallRydh2}. If $\frak X$ is N-quasi-projective over $\mathbf BG$, then $p\colon\frak X\to \mathbf BG$ is representable, so the functor $p^*\colon \QCoh(\mathbf BG)\to\QCoh(\frak X)$ preserves compact objects. Hence, locally free modules of finite rank over $\frak X$ are compact, being dualizable.
	Finally, as $\frak X$ has the resolution property, $\QCoh(\frak X)$ is generated under colimits by shifts of locally free modules of finite rank,
	by \cite[Proposition 9.3.3.7, Corollary C.2.1.7, and Corollary 9.1.3.2~(4)]{SAG}.
	Thus, we shall freely use the results of \cite{KrishnaRavi} over a general qcs base scheme $B$. 
\end{remark}

\begingroup
\def\nocontent#1#2#3{}
\let\addcontentsline=\nocontent

\subsection*{Outline}
In \S\ref{sec:MV}, we construct geometric models for classifying spaces of linear algebraic group in equivariant motivic homotopy theory. The main example is a model for the classifying space of $\GL_n$ in terms of equivariant Grassmannians.

In \S\ref{sec:periodic}, we develop some categorical machinery that will be used to equip the motivic spectrum $\KGL_\frak X$ with an $E_\infty$-ring structure. The results of this section are not otherwise essential for the proof of Theorem~\ref{thm:intro}, but they may be of independent interest.

In \S\ref{sec:KH}, we define homotopy $K$-theory of tame quotient stacks and prove that it satisfies properties (1)--(4) of Theorem~\ref{thm:intro}.

In \S\ref{sec:cdh}, we construct the motivic $E_\infty$-ring spectra $\KGL_\frak X$ representing homotopy $K$-theory and prove that they are stable under N-quasi-projective base change, which implies that $\KH$ is a cdh sheaf.

\subsection*{Notation and terminology}
This paper is a sequel to \cite{Hoyois} and uses many of the definitions and constructions introduced there, such as: the notions of homotopy invariance and Nisnevich excision \cite[Definitions 3.3 and 3.7]{Hoyois}, the corresponding localization functors $L_\htp$ and $L_\Nis$, and the combined motivic localization $L_\mot$ \cite[\S3.4]{Hoyois}; the auxiliary notion of small $G$-scheme \cite[Definition 3.1]{Hoyois}; and the definitions of the stable equivariant motivic homotopy category as a symmetric monoidal $\infty$-category and as an $\infty$-category of spectrum objects \cite[\S6.1]{Hoyois}. A notational difference with \textit{op.\ cit.} is that
we prefer to work with stacks rather than $G$-schemes, so that we write, \eg, $\SH([X/G])$ instead of $\SH^G(X)$. 

Given $\frak X\in\Tame_B$, recall that $\Sch_\frak X\subset(\Tame_B)_{/\frak X}$ is the full subcategory of N-quasi-projective $\frak X$-stacks. 
Whenever we write $\frak X$ as $[X/G]$, it is understood that $G$ is a tame affine group scheme and that $\frak X\in\Sch_{\mathbf BG}$.
If $\frak X=[X/G]$, $\Sch_\frak X$ differs slightly from the category $\Sch_X^G$ from \cite[\S3.1]{Hoyois}, but every object in either category has a Nisnevich cover whose Čech nerve belongs to their intersection, so the difference does not matter for our purposes.
We let 
$\Sm_\frak X\subset \Sch_\frak X$ be the full subcategory spanned by the smooth $\frak X$-stacks. 
We denote by $\QCoh(\frak X)^\heartsuit$ the abelian category of quasi-coherent sheaves on $\frak X$ (it is the heart of a $t$-structure on the stable $\infty$-category $\QCoh(\frak X)$ from Remark~\ref{rmk:KR}).
Given $\scr E\in\QCoh(\frak X)^\heartsuit$, we denote by $\V(\scr E)=\Spec(\Sym(\scr E))$ the associated vector bundle and by $\P(\scr E)=\Proj(\Sym(\scr E))$ the associated projective bundle.
Unless otherwise specified, presheaves and sheaves are valued in $\infty$-groupoids.

\endgroup

\section{Geometric models for equivariant classifying spaces}
\label{sec:MV}

In this section, we fix a base stack $\frak S=[S/G]\in\Tame_B$.
If $\Gamma$ is an fppf sheaf of groups on $\Sch_\frak S$, we denote by $B_\fppf\Gamma=L_\fppf(*/\Gamma)$ the presheaf of groupoids classifying $\Gamma$-torsors in the fppf topology, which we will often implicitly regard as a presheaf on $\Sm_\frak S$ (note however that the fppf sheafification must be performed on the larger category $\Sch_\frak S$).
For example, for $\frak X\in\Sm_\frak S$ and $n\geq 0$, $(B_\fppf \GL_n)(\frak X)$ is the groupoid of vector bundles of rank $n$ on $\frak X$.
When $\frak S$ is a scheme and $\Gamma$ is a smooth linear group scheme over $\frak S$, Morel and Voevodsky constructed in \cite[\S4.2]{MV} a geometric model for $L_\mot(B_\fppf\Gamma)$, \ie, they expressed $L_\mot(B_\fppf\Gamma)$ as a simple colimit of representables in $\H(\frak S)$. In this section, we generalize their result to arbitrary $\frak S\in\Tame_B$.

Let $U$ be an fppf sheaf on $\Sch_\frak S$ with an action of $\Gamma$. If $X$ is an fppf sheaf and $\pi\colon T\to X$ is a torsor under $\Gamma$, we denote by $U_\pi$ the $\pi$-twisted form of $U$, \ie,
the sheaf $L_\fppf((U\times T)/\Gamma)$. 
The Morel–Voevodsky construction is based on the following tautological lemma:

\begin{lemma}
	\label{lem:MV}
	Let $\Gamma$ be an fppf sheaf of groups on $\Sch_\frak S$ acting on an fppf sheaf $U$. Suppose that, for every $\frak X\in\Sm_\frak S$ and every fppf torsor $\pi\colon T\to \frak X$ under $\Gamma$, $U_\pi\to \frak X$ is a motivic equivalence on $\Sm_\frak S$.
	Then the map
	\[
	L_\fppf(U/\Gamma)\to B_\fppf \Gamma
	\]
	induced by $U\to *$ is a motivic equivalence on $\Sm_\frak S$.
\end{lemma}

\begin{proof}
	By universality of colimits, it suffices to show that, for every $\frak X\in\Sm_\frak S$ and every map $\frak X\to L_\fppf(*/\Gamma)$, the projection $L_\fppf(U/\Gamma\times_{*/\Gamma}\frak X)\to \frak X$ is a motivic equivalence on $\Sm_\frak S$. This is exactly the assumption.
\end{proof}

Recall from \cite[Definition 3.1]{Hoyois} that a $G$-scheme $X$ over $B$ is \emph{small} if there exists a $G$-quasi-projective morphism $X\to U$ where $U$ is affine, has trivial $G$-action, and has the $G$-resolution property. Every $\frak X\in\Tame_B$ admits a schematic Nisnevich cover whose Čech nerve consists entirely of stacks of the form $[X/G]$ where $X$ is small \cite[Lemma 3.11]{Hoyois}. 

\begin{definition}
	A \emph{system of vector bundles} over $\frak S$ is a diagram of vector bundles $(V_i)_{i\in I}$ over $\frak S$, where $I$ is a filtered poset, whose transition maps are vector bundle inclusions. Such a system is called:
	\begin{itemize}
		\item \emph{saturated} if, for every $i\in I$, there exists $2i\geq i$ such that $V_i\into V_{2i}$ is isomorphic under $V_i$ to $(\id,0)\colon V_i\into V_i\times_\frak SV_i$.
		\item \emph{complete} if, for every $\frak X=[X/G]\in\Sch_\frak S$ with $X$ small and affine, and for every vector bundle $E$ on $\frak X$, there exists $i\in I$ and a vector bundle inclusion $E\into V_i\times_\frak S\frak X$.\footnote{This is similar to the notion of \emph{complete $G$-universe} in equivariant homotopy theory.}
	\end{itemize}
\end{definition}

Note that both properties are preserved by any base change $\frak T\to \frak S$ in $\Sch_\frak S$.
The following example shows that complete saturated systems of vector bundles always exist.

\begin{example}
	\leavevmode
	\begin{enumerate}
		\item If $G$ is finite locally free and $p\colon S\to \frak S=[S/G]$ is the quotient map, then
		$(\V(p_*\scr O_S^n))_{n\geq 0}$ is a complete saturated system of vector bundles over $\frak S$.
		\item Let $\{V_\alpha\}_{\alpha\in A}$ be a set of representatives of isomorphism classes of vector bundles over $\frak S$, let $I$ be the filtered poset of maps $A\to \mathbb N$ with finitely many nonzero values, and for $i\in I$ let $V_i=\bigoplus_{\alpha\in A}V_\alpha^{i_\alpha}$. Then $(V_i)_{i\in I}$, with the obvious transition maps, is clearly a saturated system of vector bundles over $\frak S$. It is also complete, by Lemma~\ref{lem:complete} below.
	\end{enumerate}
\end{example}

\begin{lemma}\label{lem:complete}
	Let $f\colon \frak T\to \frak S$ be a quasi-affine morphism.
	For every vector bundle $V$ on $\frak T$, there exists a vector bundle $W$ on $\frak S$ and a vector bundle inclusion $V\into W\times_\frak S\frak T$.
\end{lemma}

\begin{proof}
	Let $V=\V(\scr E)$.
	Since $f$ is quasi-affine, $f^*f_*(\scr E)\to\scr E$ is an epimorphism. Since $f_*(\scr E)$ is the union of its finitely generated quasi-coherent submodules \cite[Lemma 2.10]{Hoyois}, there exists $\scr M\subset f_*(\scr E)$ finitely generated such that $f^*(\scr M)\to \scr E$ is an epimorphism.
	By the resolution property, we may assume that $\scr M$ is locally free. Setting $W=\V(\scr M)$, we then have a vector bundle inclusion $V\into W\times_\frak S\frak T$, as desired.
\end{proof}

\begin{lemma}\label{lem:lift}
	Let $\frak X=[X/G]\in\Tame_B$ with $X$ small and affine, let $s\colon\frak Z\into\frak X$ be a closed immersion, and let $V$ be a vector bundle on $\frak X$. Then any section of $V$ over $\frak Z$ lifts to a section of $V$ over $\frak X$.
\end{lemma}

\begin{proof}
	Let $V=\V(\scr E)$. We must show that any map $\scr O_\frak X\to s_*s^*(\scr E^\vee)$ in $\QCoh(\frak X)^\heartsuit$ lifts to a map $\scr O_\frak X\to \scr E^\vee$. Since $s$ is a closed immersion, the restriction map $\scr E^\vee\to s_*s^*(\scr E^\vee)$ is an epimorphism in $\QCoh(\frak X)^\heartsuit$. Moreover, since $X$ is small and affine and $G$ is linearly reductive, $\scr O_\frak X$ is projective in $\QCoh(\frak X)^\heartsuit$ \cite[Lemma 2.17]{Hoyois}. The result follows.
\end{proof}

\begin{lemma}
	\label{lem:MV2}
	Let $(V_i)_{i\in I}$ be a saturated system of vector bundles over $\frak S$.
	For every $i\in I$, let $U_i\subset V_i$ be an open substack such that $V_i\into V_j$ maps $U_i$ to $U_j$ whenever $i\leq j$. Suppose that:
	\begin{enumerate}
		\item 
		there exists $i\in I$ such that $U_{i}\to\frak S$ has a section;
		\item for all $i\in I$, under the isomorphism $V_{2i}\simeq V_i^2$, $(U_i\times V_i)\cup (V_i\times U_i)\subset U_{2i}$. 
	\end{enumerate}
	Then $U_\infty= \colim_{i\in I}U_i\in\PSh(\Sm_\frak S)$ is motivically contractible.
\end{lemma}

\begin{proof}
	By \cite[Proposition 3.16 (2)]{Hoyois}, it will suffice to show that, for every $\frak X=[X/G]\in\Sm_\frak S$ with $X$ small and affine, the simplicial set $\Map(\A^\bullet\times\frak X,U_\infty)$ is a contractible Kan complex. Consider a lifting problem
	\begin{tikzmath}
		\diagram{\del\Delta^n & \Map(\A^\bullet\times\frak X,U_\infty)\rlap. \\ \Delta^n & \\};
		\arrows (11-) edge node[above]{$f$} (-12) (11) edge[c->] (21) (21) edge[dashed] (12);
	\end{tikzmath}
	Then $f\colon\del\A^n_\frak X\to U_\infty$ is a morphism from the boundary of the algebraic $n$-simplex over $\frak X$ to $U_\infty$, and it factors through $U_i$ for some $i$ since $\del\A^n_\frak X$ is compact as an object of $\PSh(\Sm_\frak S)$. Increasing $i$ if necessary, we may assume, by (1), that there exists a section $x\colon\frak S\to U_i$.
	By Lemma~\ref{lem:lift}, there exists a morphism $g\colon \A^n_\frak X\to V_i$ lifting $f$.
	Choose a closed substack $Z_i\subset V_i$ complementary to $U_i$, so that $g^{-1}(Z_i)\cap\del\A^n_\frak X=\emptyset$. Again by Lemma~\ref{lem:lift}, the map 
	\[
	g^{-1}(Z_i)\sqcup \del\A^n_\frak X \to \frak S\sqcup\frak S \xrightarrow{x\sqcup 0} V_i
	\]
	admits an extension $h\colon \A^n_\frak X\to V_i$. Then the morphism $(g,h)\colon \A^n_\frak X\to V_{i}^2$ misses $Z_i^2$ and hence solves the lifting problem, by (2).
\end{proof}

If $E$ is a vector bundle over $\frak S$, we denote by $\GL(E)$ the group $\frak S$-stack of linear automorphisms of $E$. By a \emph{subgroup} of $\GL(E)$ we mean a subfunctor of its functor of points (valued in group objects in groupoids).

\begin{theorem}
	\label{thm:MV}
	Let $E$ be a vector bundle over $\frak S$, $\Delta\subset\GL(E)$ a closed subgroup, and $\Gamma\subset\Delta$ a subgroup that is flat and finitely presented over $\frak S$.
	Let $(V_i)_{i\in I}$ be a complete saturated system of vector bundles over $\frak S$. For each $i\in I$, let $U_i\subset \Hom(E,V_i)$ be the open substack where the action of $\Delta$ is strictly free,
	and let $U_\infty=\colim_{i\in I} U_i$. Then the map
	\[
	L_\fppf(U_\infty/\Gamma)\to B_\fppf \Gamma
	\]
	induced by $U_\infty\to *$ is a motivic equivalence on $\Sm_\frak S$. 
\end{theorem}

\begin{proof}
	We check that $U_\infty$ satisfies the assumption of Lemma~\ref{lem:MV}, \ie, that for any $\frak X\in\Sm_\frak S$ and any $\Gamma$-torsor $\pi\colon T\to\frak X$, the map $(U_\infty)_\pi\to\frak X$ is a motivic equivalence on $\Sm_\frak S$. By \cite[Proposition 4.6]{Hoyois}, we can assume that $\frak X=[X/G]$ with $X$ small and affine.
	 It then suffices to show that the saturated system of vector bundles $\Hom(E_\pi,V_i\times_\frak S\frak X)$ over $\frak X$ and the open substacks $(U_i)_\pi$ satisfy the conditions of Lemma~\ref{lem:MV2} with $\frak S=\frak X$.
	The second condition is clear, by definition of $U_i$. To verify the first condition, we can assume that $\Delta=\GL(E)$.
	Sections of $(U_i)_\pi$ over $\frak X$ are then vector bundle inclusions $E_\pi\into V_i\times_\frak S\frak X$. Since $(V_i)_{i\in I}$ is complete, there exist such inclusions for large enough $i$.
\end{proof}

\begin{remark}
	Although this is not always true in the generality of Theorem~\ref{thm:MV}, the fppf quotients $L_\fppf(U_i/\Gamma)$ are often representable by (necessarily smooth) quasi-projective $\frak S$-stacks, so that the presheaf $L_\fppf(U_\infty/\Gamma)$ is a filtered colimit of representables. It is in that sense that $L_\fppf(U_\infty/\Gamma)$ is a geometric model for $L_\mot (B_\fppf\Gamma)$.
\end{remark}

\begin{corollary}\label{cor:BC}
	Under the assumptions of Theorem~\ref{thm:MV}, suppose that the fppf quotients $L_\fppf(U_i/\Gamma)$ are universally representable by N-quasi-projective $\frak S$-stacks. Then, for every N-quasi-projective morphism $f\colon\frak T\to\frak S$, the map
	 \[f^*(B_\fppf\Gamma)\to B_\fppf(f^*\Gamma)\]
	 in $\PSh(\Sm_\frak T)$ is a motivic equivalence.
\end{corollary}

\begin{proof}
	Consider the following commutative square in $\PSh(\Sm_\frak T)$:
	\begin{tikzmath}
		\diagram{f^*L_\fppf(U_\infty/\Gamma) & f^*(B_\fppf\Gamma) \\
		L_\fppf(f^*(U_{\infty}/\Gamma)) & B_\fppf(f^*\Gamma)\rlap. \\};
		\arrows (11-) edge (-12) (11) edge (21) (21-) edge (-22) (12) edge (22);
	\end{tikzmath}
	By Theorem~\ref{thm:MV}, the horizontal maps are motivic equivalences. On the other hand, by assumption, the left vertical arrow is an isomorphism between ind-representable presheaves on $\Sm_\frak T$.
\end{proof}

Corollary~\ref{cor:BC}, applied to $\Gamma=\GL_n$, is all that we will need from this section in the sequel. In that case, $U_i\subset\Hom(\A^n_\frak S,V_i)$ is the open substack of vector bundle inclusions, and $L_\fppf(U_i/\GL_n)$ is universally represented by the Grassmannian $\Gr_n(V_i)$. Let us make Theorem~\ref{thm:MV} more explicit in this special case:

\begin{corollary}\label{cor:BGLn}
	Let $(V_i)_{i\in I}$ be a complete saturated system of vector bundles over $\frak S$.
	For any $n\geq 0$, the map
	\[
	\colim_{i\in I} \Gr_n(V_i) \to B_\fppf\GL_n
	\]
	in $\PSh(\Sm_\frak S)$ classifying the tautological bundles is a motivic equivalence.
\end{corollary}

\section{Periodic $E_\infty$-algebras}
\label{sec:periodic}

Let $\scr C\in \CAlg(\Pr^\mathrm{L})$ be a presentably symmetric monoidal $\infty$-category, $S$ a set of objects of $\scr C_{/\1}$, and $\scr M$ a $\scr C$-module in $\Pr^\mathrm{L}$. For every $x\in\scr C$, we have the adjunction
\[
x\otimes\ph: \scr M \rightleftarrows \scr M: \Hom(x,\ph).
\]
We say that an object $E\in\scr M$ is \emph{$S$-periodic} if $\alpha^*\colon E\to \Hom(x,E)$ is an equivalence for every $\alpha\colon x\to\1$ in $S$. We denote by $P_S\scr M\subset \scr M$ the full subcategory spanned by the $S$-periodic objects.
It is clear that this inclusion is an accessible localization and hence admits a left adjoint $P_S$, called \emph{periodization}.
Note that $E$ is $S$-periodic if and only if it is local with respect to $\id_M\tens\alpha$ for every $M\in\scr M$ and $\alpha\in S$. If $\scr M=\scr C$, it follows immediately that the localization functor $P_S$ is compatible with the monoidal structure in the sense of \cite[Definition 2.2.1.6]{HA}, and hence that it can be promoted to a symmetric monoidal functor \cite[Proposition 2.2.1.9]{HA}. In particular, for every $E_\infty$-algebra $A$ in $\scr C$, $P_SA$ is also an $E_\infty$-algebra in $\scr C$ and $A\to P_SA$ is an $E_\infty$-map.

Let $S_0$ be the set of domains of morphisms in $S$.
 Consider the presentably symmetric monoidal $\infty$-category $\scr C[S_0^{-1}]$ obtained from $\scr C$ by adjoining formal inverses to elements of $S_0$ (see \cite[\S6.1]{Hoyois}), which is in particular a $\scr C$-module.
We have an adjunction
\begin{tikzmath}
	\diagram{\scr C & \scr C[S_0^{-1}]\rlap, \\};
	\arrows (11-) edge[vshift=\dbl] node[above]{$\Phi$} (-12) (-12) edge[vshift=\dbl] node[below]{$\Psi$} (11-);
\end{tikzmath}
where $\Phi$ is symmetric monoidal. It follows that $\Psi$ preserves $S$-periodic objects. Hence, the above adjunction induces an adjunction
\begin{tikzequation}\label{eqn:Pdeloop}
	\diagram{P_S\scr C & P_S(\scr C[S_0^{-1}])\rlap. \\};
	\arrows (11-) edge[vshift=\dbl] node[above]{$P_S\Phi$} (-12) (-12) edge[vshift=\dbl] node[below]{$\Psi$} (11-);
\end{tikzequation}

\begin{proposition}\label{prop:Pdeloop}
	Let $\scr C$ be a presentably symmetric monoidal $\infty$-category, $S$ a set of objects of $\scr C_{/\1}$, and $S_0$ the set of domains of morphisms in $S$. Then the adjunction~\eqref{eqn:Pdeloop} is an equivalence of symmetric monoidal $\infty$-categories.
	In particular, every $S$-periodic $E_\infty$-algebra in $\scr C$ lifts uniquely to an $S$-periodic $E_\infty$-algebra in $\scr C[S_0^{-1}]$.
\end{proposition}

\begin{proof}
	Indeed, the symmetric monoidal functors $P_S\colon\scr C\to P_S\scr C$ and $P_S\Phi\colon\scr C\to P_S(\scr C[S_0^{-1}])$ satisfy the same universal property, since the former sends every $x\in S_0$ to an invertible object, namely, the unit of $P_S\scr C$.
\end{proof}

We would like to understand the periodization functor $P_S$ more explicitly.
Consider the case where $S$ consists of a single map $\alpha\colon x\to\1$. Given $E\in\scr C$, it is tempting to think that $P_\alpha E$ is given by the formula
\[
\colim \bigl( E \stackrel\alpha\longrightarrow \Hom(x,E) \stackrel \alpha\longrightarrow \Hom(x^{\tens 2},E) \stackrel \alpha\longrightarrow \dotsb\bigr),
\]
at least if we assume that $\Hom(x,\ph)$ preserves filtered colimits (otherwise, we would naturally consider a transfinite construction). This formula is indeed correct if $\scr C$ is a stable $\infty$-category and $\alpha\colon\1\to\1$ is multiplication by an integer, but not in general. For example, suppose that $\scr C$ is the symmetric monoidal $\infty$-category of small stable $\infty$-categories, and let $\alpha$ be multiplication by a positive integer on the unit $\Sp^\fin$. Then $P_\alpha\scr C\subset\scr C$ is the subcategory of zero objects, but the above colimit with $E=\Sp^\fin$ is not zero. The essential difference between these two cases is the following: in the first case, the cylic permutation of $\alpha^3$ is homotopic to the identity (because it is the image of an even element in $\pi_1$ of the sphere spectrum), but in the second case, no nontrivial permutation of $\alpha^n$ is homotopic to the identity.
We will show that there exists an analogous formula for $P_S$ in general, provided that the elements of $S$ are cyclically symmetric in a suitable sense.

We recall some constructions from \cite[\S6.1]{Hoyois}.
Let $X$ be any set of objects of $\scr C$.
The filtered simplicial set $L(X)$ is the union over finite subsets $F\subset X$ of the simplicial sets $L^F$, where $L$ is the $1$-skeleton of the nerve of the poset $\N$. We view a vertex of $L(X)$ as a formal tensor product of elements of $X$. The $\scr C$-module $\Stab_{X}(\scr C)$ of \emph{$X$-spectra} is then defined as the limit of a diagram $L(X)^\op\to \Mod_\scr C$ taking each vertex of $L(X)$ to $\scr C$ and each arrow $w\to w\tens x$ to the functor $\Hom(x,\ph)$. Equivalently, $\Stab_X(\scr C)$ is the $\infty$-category of cartesian sections of the cartesian fibration over $L(X)$ classified by $L(X)^\op\to\Cat_\infty$. A general section of this cartesian fibration will be called an \emph{$X$-prespectrum} in $\scr C$; we denote by $\Stab_{X}^\lax(\scr C)$ the $\scr C$-module of $X$-prespectra. Thus, $\Stab_{X}(\scr C)$ is a (left exact) localization of $\Stab_X^\lax(\scr C)$. The localization functor is called \emph{spectrification} and is denoted by $Q\colon \Stab^\lax_X(\scr C)\to \Stab_X(\scr C)$. If $\Hom(x,\ph)$ preserves filtered colimits for all $x\in X$, which will be the case in all our applications, spectrification is given by the familiar formula
\[
Q(E)_w = \colim_{v\in L(X)} \Hom(v, E_{w\tens v}).
\]
In general, one can describe spectrification as follows. For every $x\in X$, consider the full subcategory $\scr E_x\subset \Stab_X^\lax(\scr C)$ consisting of $X$-prespectra that are spectra in the $x$-direction, so that $\Stab_X(\scr C)=\bigcap_{x\in X}\scr E_x$.
Choose a regular cardinal $\kappa$ such that $\Hom(x,\ph)$ preserves $\kappa$-filtered colimits for all $x\in X$, and let $\sh_x$ be the pointed endofunctor of $\Stab_X^\lax(\scr C)$ given by $\sh_x(E)_w=\Hom(x,E_{w\tens x})$.
Then the $\kappa$th iteration $\sh_x^{\kappa}$ of $\sh_x$ lands in $\scr E_x$.
Moreover, any map $E\to F$ with $F\in\scr E_x$ factors uniquely through $\sh_x(E)$. It follows that $\sh_x^{\kappa}$ is left adjoint to the inclusion $\scr E_x\subset \Stab_X^\lax(\scr C)$. The total localization functor $Q$ can now be written as an appropriate $\kappa$-filtered transfinite composition in which each indecomposable map is an instance of $\id\to \sh_x^{\kappa}$ for some $x\in X$ (see the proof of \cite[Lemma 7.3.2.3]{HTT}).

To every $E\in\scr C$ we can associate a ``constant'' $S_0$-prespectrum $c_{S}E=(E)_{w\in L(S_0)}$ with structure maps $E\to\Hom(x,E)$ induced by the maps in $S$. 
Let $Q_S\colon\scr C \to \scr C$ be the functor defined by 
\[
Q_SE=\Omega^\infty Q(c_SE),
\]
where $\Omega^\infty\colon \Stab_{S_0}(\scr C)\to\scr C$ is evaluation at the initial vertex of $L(S_0)$.
There is an obvious natural transformation $\id\to Q_S$. For example, if $S$ consists of a single map $\alpha\colon x\to\1$ and $\Hom(x,\ph)$ preserves filtered colimits, we have
\[
Q_\alpha E = \colim\bigl(E \xrightarrow{\alpha} \Hom(x,E) \xrightarrow\alpha \Hom(x^{\tens 2}, E)\to\dotsb\bigr).
\]

\begin{lemma}\label{lem:P=Q}
	Let $\scr C$ be a presentably symmetric monoidal $\infty$-category, $S$ a set of objects of $\scr C_{/\1}$, and $E\in\scr C$. If $Q_SE$ is $S$-periodic, then the map $E\to Q_SE$ exhibits $Q_SE$ as the $S$-periodization of $E$.
\end{lemma}

\begin{proof}
	For $x\in S_0$, the functor $\Hom(x,\ph)\colon P_S\scr C\to P_S\scr C$ is an equivalence of $\infty$-categories, since $P_S(x)$ is invertible in $P_S\scr C$. Hence, $\Omega^\infty\colon P_S\Stab_{S_0}(\scr C)\to P_S\scr C$ is an equivalence.
		Consider the following commutative diagram of $\scr C$-modules:
		\begin{tikzmath}
			\diagram{\scr C & \Stab^\lax_{S_0}(\scr C) & \Stab_{S_0}(\scr C) & \\
			P_S\scr C & P_S\Stab^\lax_{S_0}(\scr C) &  P_S\Stab_{S_0}(\scr C) & P_S\scr C\rlap.\\};
			\arrows (11-) edge node[below]{$\Sigma^\infty_\lax$} (-12) (12-) edge node[below]{$Q$} (-13)
			(21-) edge (-22) (22-) edge (-23) (23-) edge node[above]{$\Omega^\infty$} node[below]{$\simeq$} (-24)
			(11) edge (21) (12) edge (22) (13) edge (23);
			\draw[->,font=\scriptsize] (11) to[out=15,in=165] node[above]{$\Sigma^\infty$} (13);
			\draw[->,font=\scriptsize] (21) to[out=-15,in=-165] node[below]{$\simeq$} ([yshift=4pt]23.south west);
		\end{tikzmath}
		All the vertical arrows are periodization functors, and the lower composition is the identity. This diagram shows that
		\[
		P_S(E) = \Omega^\infty P_SQ(\Sigma^\infty_\lax E).
		\]
		Here, $\Sigma^\infty_\lax E$ is the free $S_0$-prespectrum $(E\tens w)_{w\in L(S_0)}$.
		The obvious map $\Sigma^\infty_\lax E\to c_SE$ is manifestly a termwise $P_S$-equivalence.
		Since the right adjoints to the various evaluation functors $\Stab_{S_0}^\lax(\scr C)\to\scr C$ preserve $S$-periodic objects, termwise $P_S$-equivalences of $S_0$-prespectra are in fact $P_S$-equivalences. It follows that
		\[
		P_S(E) = \Omega^\infty P_SQ(c_S E).
		\]
		All the terms of the $S_0$-spectrum $Q(c_SE)$ are equivalent to $Q_SE$. Hence, by the assumption, $Q(c_SE)$ is already $S$-periodic, and we get $P_SE=Q_SE$, as desired.
\end{proof}

\begin{example}\label{ex:KB}
	Let $K$ denote the presheaf of $E_\infty$-ring spectra $X\mapsto K(X)$ on qcqs schemes, and let $\beta\in\tilde K_1(\G_m,1)$ be the \emph{Bott element}, that is, the element induced by the automorphism $t$ of $\scr O_{\G_m}$, where $\G_m=\Spec \Z[t^{\pm 1}]$.
	Let $\gamma$ be the composite
	\[
	(\P^1\minus 0)\coprod_{\G_m}\A^1 \to \Sigma(\G_m/1) \xrightarrow{\beta} K,
	\]
	where the pushout is taken in presheaves and pointed at $1$.
	By inspecting the definition \cite[Definition 6.4]{TT}, we see that the Bass–Thomason–Trobaugh $K$-theory spectrum $K^B$ is the $K$-module $Q_\gamma K$.
	Since $K^B$ is $\gamma$-periodic, Lemma~\ref{lem:P=Q} implies that $K^B=P_\gamma K$.
	In particular, $K^B$ is an $E_\infty$-algebra under $K$. The same argument applies to $K$ and $K^B$ as presheaves on $\Tame_B$ (see \cite[\S3.5]{KrishnaRavi} for the definition of $K^B$ in this context).
\end{example}

Let $\scr C$ be a symmetric monoidal $\infty$-category.
An object $x\in \scr C$ is called \emph{$n$-symmetric} if the cyclic permutation $\sigma_n$ of $x^{\tens n}$ is homotopic to the identity. We will say that $x$ is \emph{symmetric} if it is $n$-symmetric for some $n\geq 2$.
If $\scr C$ is presentably symmetric monoidal and $X$ is a set of symmetric objects of $\scr C$, there is an equivalence of $\scr C$-modules $\scr C[X^{-1}]\simeq \Stab_{X}(\scr C)$ (see \cite[Corollary 2.22]{Robalo} and \cite[\S6.1]{Hoyois}).
 
The $\infty$-category $\scr C_{/\1}$ inherits a symmetric monoidal structure from $\scr C$ such that the forgetful functor $\scr C_{/\1}\to\scr C$ is symmetric monoidal. An $n$-symmetric object in $\scr C_{/\1}$ is then a morphism $\alpha\colon x\to\1$ such that the cyclic permutation $\sigma_n$ of $x^{\tens n}$ is homotopic \emph{over $\1$} to the identity. 

\begin{example}
	If $\scr C$ is symmetric monoidal, $\End(\1)$ is an $E_\infty$-space under composition. In particular, for every $\alpha\colon \1\to\1$, the cyclic permutation of $n$ letters induces a self-homotopy $\sigma_n$ of $\alpha^n$.
	Then $\alpha$ is $n$-symmetric in $\scr C_{/\1}$ if and only if the homotopy class of $\sigma_n$ is in the image of the group homomorphism $\pi_1(\End(\1),\id)\to \pi_1(\End(\1),\alpha^n)$ induced by $\End(\1)\to \End(\1)$, $\beta\mapsto \alpha^n\circ\beta$. In particular, if $\sigma_n$ vanishes in $\pi_1(\End(\1),\alpha^n)$, then $\alpha$ is $n$-symmetric in $\scr C_{/\1}$.
\end{example}

\begin{lemma}\label{lem:symmetry}
	Let $\scr C$ be a symmetric monoidal $\infty$-category and let $\alpha\colon x\to\1$ be a symmetric object in $\scr C_{/\1}$. Let $X_\bullet$ be the tower $\N^\op\to \scr C_{/\1}$, $k\mapsto x^{\tens k}$, with transition maps $\id\tens\alpha$. Then the transformations \[\alpha\tens \id,\;\id\tens\alpha\colon X_{\bullet+1}\to X_\bullet\] are homotopic as maps in $\mathrm{Pro}(\scr C_{/\1})$. 
\end{lemma}

\begin{proof}
	Let $\sigma_k$ be the cyclic permutation of $x^{\tens k}$ that moves the first factor to the end.
	The map $\id\tens \alpha\colon x^{\tens k+1} \to x^{\tens k}$ is then the composite of $\sigma_{k+1}^{-1}$ and $\alpha\tens\id$. Define a new tower $\tilde X_\bullet\colon \N^\op\to\scr C$
	with $\tilde X_k=x^{\tens k}$ and with transition maps $\sigma_{k}^{-1}\circ (\id\tens\alpha)\circ\sigma_{k+1}\colon x^{\tens k+1}\to x^{\tens k}$. The permutations $\sigma_k$ assemble into a natural equivalence $\sigma\colon\tilde X_\bullet \to X_\bullet$ such that $(\id\tens\alpha)\circ\sigma\simeq\alpha\tens\id$.
	The strategy of the proof is the following: we will construct an equivalence of pro-objects $\zeta\colon X_{\bullet+1}\to \tilde X_{\bullet+1}$ making the diagram
	\begin{tikzmath}
		\def\rowsep{4em}
		\diagram{X_{\bullet + 1} & \tilde X_{\bullet+1} & X_{\bullet+1} \\
		& X_\bullet & \\};
		\arrows (11-) edge node[above]{$\zeta$} (-12) (12-) edge node[above]{$\sigma$} (-13)
		(11) edge node[below left]{$\alpha\tens\id$} (22) (12) edge node[fill=white,pos=.4]{$\alpha\tens\id$} (22) (13) edge node[below right]{$\id\tens\alpha$} (22);
	\end{tikzmath}
	commute and such that $\sigma\circ\zeta$ is homotopic to the identity.
	Let us call $\pi$ and $\tilde\pi$ the morphisms $\alpha\tens\id\colon X_{\bullet+1}\to X_\bullet$ and $\alpha\tens\id\colon \tilde X_{\bullet+1}\to X_\bullet$.
	
	Suppose that $\alpha$ is $(n+1)$-symmetric, and let $L$ be the $1$-skeleton of the nerve of the poset $n\N\subset \N$.
	We will then construct $\zeta$ as a morphism in $\Fun(L^\op,\scr C)$, and we will prove that $\tilde\pi\circ\zeta\simeq\pi$ and $\sigma\circ\zeta\simeq\id$ in $\Fun(L^\op,\scr C)$.
	The image of an edge of $L$ by either $\pi$ or $\tilde\pi$ has the form
	\begin{tikzmath}
		\def\colsep{5em}
		\diagram{
		x^{\tens nk + 1} & x^{\tens n(k-1)+1} \\
		x^{\tens nk}     & x^{\tens n(k-1)}\rlap, \\
		};
		\arrows (11-) edge (-12) (11) edge node[left]{$\alpha\tens\id$} (21) (21-) edge node[below]{$\id\tens\alpha^{n}$} (-22) (12) edge node[right]{$\alpha\tens\id$} (22)
		(11) edge node[fill=white]{$\alpha\tens\id\tens\alpha^{n}$} (22);
	\end{tikzmath}
	but $\pi$ and $\tilde\pi$ differ on the upper triangle.
	Let $\sigma'_{k}\colon x^{\tens nk+1}\to x^{\tens nk+1}$ be the cyclic permutation $\sigma_{n+1}$ applied to the $n+1$ factors of $x^{\tens nk+1}$ that are killed by the diagonal. Observe that
	\begin{equation}\label{eqn:sigmas}
	\sigma_{nk+1} = (\sigma_{n(k-1)+1}\tens\id) \circ \sigma'_{k}.
	\end{equation}
	In particular, the transition map $x^{\tens nk + 1} \to x^{\tens n(k-1)+1}$ in $\tilde X_{\bullet+1}$ is $(\id\tens\alpha^n)\circ\sigma'_{k}$. We define $\zeta\colon X_{\bullet+1}\to\tilde X_{\bullet +1}$ to be the identity on each vertex of $L$ and the given homotopy $\sigma'\simeq\id$ on each edge. Thus, the image by $\zeta$ of an edge of $L$ is the square
	\begin{tikzmath}
		\diagram[row sep={2em,between origins},column sep={2em,between origins},nodes={inner sep=1pt,circle}]{
		\cdot & & & \cdot \\ & & & \\ \cdot & \cdot & & \cdot \\
		};
		\arrows
		(11-) edge node[above]{$\id\tens\alpha^n$} (-14) (11) edge[-] (31) (11) edge[-] (32)
		(14) edge[-] (34)
		(31-) edge node[below]{$\sigma'$} (-32)
		(32-) edge node[below]{$\id\tens\alpha^n$} (-34)
		;
	\end{tikzmath}
	where untipped lines represent identity morphisms and the triangle is the given homotopy $\sigma'\simeq\id$.
	The composites $\tilde\pi\circ\zeta$ and $\sigma\circ\zeta$ are then described by the following pictures:
	\[
	\begin{tikzpicture}
		\diagram[row sep={2em,between origins},column sep={2em,between origins},nodes={inner sep=1pt,circle}]{
		\cdot & & \cdot \\ & & \\ \cdot & \cdot & \cdot \\ & & \\ \cdot & & \cdot \\
		};
		\arrows
		(11-) edge node[above]{$\id\tens\alpha^n$} (-13) (11) edge[-] (31) (11) edge[-] (32)
		(13) edge[-] (33)
		(31-) edge node[above]{$\sigma'$} (-32) (31) edge node[left]{$\alpha\tens\id$} (51) (31) edge[shorten >=1.5pt] (53)
		(32-) edge (-33) (32) edge[shorten >=5pt] (53)
		(33) edge node[right]{$\alpha\tens\id$} (53)
		(51-) edge node[below]{$\id\tens\alpha^n$} (-53)
		;
	\end{tikzpicture}
	\qquad\qquad
	\begin{tikzpicture}
		\diagram[row sep={2em,between origins},column sep={2em,between origins},nodes={inner sep=1pt,circle}]{
		\cdot & & \cdot \\ & & \\ \cdot & \cdot & \cdot \\ \cdot & & \\ \cdot & & \cdot \\
		};
		\arrows
		(11-) edge node[above]{$\id\tens\alpha^n$} (-13) (11) edge[-] (31) (11) edge[-] (32)
		(13) edge[-] (33)
		(31-) edge node[above]{$\sigma'$} (-32) (31) edge node[left]{$\sigma'$} (41)
		(32-) edge (-33) (32) edge[-] (41)
		(33) edge node[right]{$\sigma$} (53)
		(41) edge node[left]{$\sigma\tens\id$} (51)
		(51-) edge node[below]{$\id\tens\alpha^n$} (-53)
		;
	\end{tikzpicture}
	\]
	In the first picture, the two diagonal arrows are $\alpha\tens\id\tens\alpha^n$.
	The assumption that the given homotopy $\sigma_{n+1}\simeq\id$ is a homotopy over $\1$ implies that the triangle with median $\sigma'$ is homotopic rel its boundary to an identity $2$-cell, showing that $\tilde\pi\circ\zeta\simeq\pi$.
	
	Using~\eqref{eqn:sigmas}, we inductively construct homotopies $\sigma_{nk+1}\simeq\id$ for $k\geq 0$.
	The pentagon in the second picture is the tensor product
	\[
	\delimiterfactor=0
	\delimitershortfall=20pt
	\left(
	\begin{tikzpicture}
		\diagram[row sep={4em,between origins}]{x^{\tens n(k-1)+1} & x^{\tens n(k-1)+1} \\ x^{\tens n(k-1)+1} & x^{\tens n(k-1)+1} \\};
		\arrows (11-) edge node[above]{$\id$} (-12) (11) edge node[left]{$\sigma_{n(k-1)+1}$} (21) (21-) edge node[below]{$\id$} (-22) (12) edge node[right]{$\sigma_{n(k-1)+1}$} (22);
	\end{tikzpicture}
	\right)
	\tens
	\left(
	\begin{tikzpicture}
		\diagram[row sep={4em,between origins}]{x^{\tens n} & \1 \\ x^{\tens n} & \1 \\};
		\arrows (11-) edge node[above]{$\alpha^n$} (-12) (11) edge node[left]{$\id$} (21) (21-) edge node[below]{$\alpha^n$} (-22) (12) edge node[right]{$\id$} (22);
	\end{tikzpicture}
	\right).
	\]
	Using the homotopies $\sigma_{n(k-1)+1}\simeq\id$ and $\sigma_k'\simeq \id$, we obtain for every edge $e\colon\Delta^1\to L$ a homotopy in $\Fun(\Delta^1,\scr C)$ between $(\sigma\circ\zeta)_e$ and the identity. By construction, these homotopies agree on the common vertex of two consecutive edges of $L$, and hence they define a homotopy $\sigma\circ\zeta\simeq\id$, as desired.
\end{proof}

\begin{theorem}\label{thm:P=Q}
	Let $\scr C$ be a presentably symmetric monoidal $\infty$-category and $S$ a set of symmetric objects of $\scr C_{/\1}$. Then $P_S\simeq Q_S$. More precisely, for every $E\in\scr C$, the canonical map $E\to Q_SE$ exhibits $Q_SE$ as the $S$-periodization of $E$.
\end{theorem}

\begin{proof}
	By Lemma~\ref{lem:P=Q}, it suffices to show that $Q_SE$ is $S$-periodic, \ie, that $c_SQ_SE$ is an $S_0$-spectrum.
	We use the following explicit description of the spectrification functor $Q$, from the proof of \cite[Lemma 7.3.2.3]{HTT}.
	Choose a regular cardinal $\kappa$ such that $\Hom(x,\ph)$ preserves $\kappa$-filtered colimits for all $x\in S_0$, and choose a bijection $f\colon S_0\to\lambda$ for some ordinal $\lambda$. Then $Q=\colim_{\mu<\lambda\kappa} F_\mu$, where $F_{\mu+1}=\sh_x^\kappa F_\mu$ if $\mu=\lambda\nu+f(x)$.
	
	Note that each $S_0$-prespectrum $F_\mu (c_SE)$ is ``constant'' in the sense that all its terms and structure maps in a given direction are the same.
	For any $\alpha\colon x\to\1$ and $\beta\colon y\to\1$ in $S$ with $\alpha\neq\beta$, it is clear that the structure map of $\sh_x(c_SE)$ in the $y$-direction is $\beta^*$.
	Lemma~\ref{lem:symmetry} shows that the structure map of $\sh_x^\omega (c_SE)$ in the $x$-direction is naturally homotopic to $\alpha^*$ under $E$. Hence, we have an equivalence $\sh_x^\omega (c_SE)\simeq c_S\Omega^\infty \sh_x^\omega (c_SE)$ under $c_SE$. By a straightforward transfinite induction, we can identify the towers $\{F_\mu (c_SE)\}_{\mu\leq\lambda\kappa}$ and $\{c_S\Omega^\infty F_\mu (c_SE)\}_{\mu\leq\lambda\kappa}$. In particular, $Q(c_SE)\simeq c_SQ_SE$ and $c_SQ_SE$ is an $S_0$-spectrum.
\end{proof}

\section{Homotopy $K$-theory of tame quotient stacks}
\label{sec:KH}

The homotopy $K$-theory spectrum $\KH(X)$ of a qcqs scheme $X$ is the geometric realization of the simplicial spectrum $K^B(\Delta^\bullet\times X)$, where $K^B$ is the Bass–Thomason–Trobaugh $K$-theory. Equivalently,
\[
\KH=L_{\A^1}K^B,
\]
where $L_{\A^1}$ is the reflection onto the subcategory of $\A^1$-homotopy invariant presheaves (often called the naive $\A^1$-localization).
There is an alternative point of view on $\KH$ due to Cisinski \cite{Cisinski}. 
An important feature of the Bass construction is that $K^B$ is a Nisnevich sheaf, whereas $K$ is not. It is also clear that the naive $\A^1$-localization functor $L_{\A^1}$ preserves Nisnevich sheaves of spectra, so that $\KH$ is not only $\A^1$-invariant but is also a Nisnevich sheaf. The canonical map $K\to \KH$ therefore factors through the so-called motivic localization $L_\mot(K)=L_{\A^1}L_\Nis(K)$. But the resulting map $L_\mot(K)\to\KH$ is not yet an equivalence: instead, it exhibits $\KH$ as the periodization of $L_\mot(K)$ with respect to the Bott element $\beta\in \tilde K_1(\G_m,1)$.
Our definition of the homotopy $K$-theory of a stack $\frak X$ is directly analogous to this construction. The main difficulty is that we now have to deal with several Bott elements: one for each vector bundle over $\frak X$. We also have to replace $L_{\A^1}$ by the more complicated homotopy localization $L_\htp$ \cite[\S3.2]{Hoyois}, which, unlike $L_{\A^1}$, need not preserve Nisnevich sheaves of spectra. Nevertheless, we will see that the identity $\KH=L_{\A^1}K^B$ still holds for quotient stacks $[X/G]$ with $G$ finite or diagonalizable. 

For $\frak X\in\Tame_B$, 
we will denote by $K_\frak X$ and $K_\frak X^B$ the restrictions of $K$ and $K^B$ to $\Sch_\frak X$.	
Let $\scr E$ be a locally free module of finite rank $r$ over $\frak X$,\footnote{We do not assume that $\scr E$ has constant rank, so $r$ is a locally constant integer on $\frak X$. Formulas involving $r$ must be interpreted accordingly.}
 $\P(\scr E)$ the associated projective bundle, and $\scr O(1)$ the universal sheaf on $\P(\scr E)$.
 By the projective bundle formula \cite[Theorem 3.6]{KrishnaRavi}, the functors\footnote{Here, $\Perf(\frak Y)$ is the stable $\infty$-category of perfect complexes over $\frak Y$, \ie, dualizable objects in $\QCoh(\frak Y)$.}
 \begin{equation}\label{eqn:Perf}
 \Perf(\frak Y)\to \Perf(\frak Y\times_\frak X\P(\scr E)), \quad E\mapsto E\boxtimes\scr O(-i),
 \end{equation}
 for $\frak Y\in\Sch_\frak X$ and $0\leq i\leq r-1$,
 induce an equivalence of $K_\frak X$-modules
\[
\prod_{i=0}^{r-1} K_\frak X \simeq \Hom(\P(\scr E)_+,K_\frak X).
\]
Let $\V^+(\scr E)$ denote the quotient $\P(\scr E\oplus\scr O_\frak X)/\P(\scr E)$, viewed as a pointed presheaf on $\Sch_\frak X$.
The right square in the following diagram is then commutative, and we get an equivalence as indicated:
\begin{tikzequation}\label{eqn:PBF}
	\diagram{\Hom(\V^+(\scr E),K_\frak X) & \Hom(\P(\scr E\oplus\scr O_\frak X)_+,K_\frak X) & \Hom(\P(\scr E)_+, K_\frak X) \\
	K_\frak X & \prod_{i=0}^{r} K_\frak X & \prod_{i=0}^{r-1}K_\frak X\rlap.\\};
	\arrows (11-) edge (-12) (12-) edge (-13)
	(11) edge[dashed] node[left]{$\simeq$} (21) (12) edge[<-] node[left]{$\simeq$} (22) (13) edge[<-] node[left]{$\simeq$} (23)
	(21-) edge[<-] node[below]{$\mathrm{pr}_r$} (-22) (22-) edge[<-] (-23);
\end{tikzequation}
This equivalence is a morphism of $K_\frak X$-modules and is therefore determined by a map $\beta_{\scr E}\colon \V^+(\scr E)\to K_\frak X$.
A standard representative of $\beta_\scr E$ by a perfect complex is given by the Koszul complex of the composition
\begin{equation}\label{eqn:Koszul}
\scr E_{\P(\scr E\oplus \scr O)}(-1) \into (\scr E\oplus \scr O)_{\P(\scr E\oplus \scr O)}(-1) \onto \scr O_{\P(\scr E\oplus \scr O)},
\end{equation}
where the first map is the inclusion of the first summand and the second map is the tautological epimorphism, tensored with $\det(\scr E)[r]^\vee$. In particular, the image of $\beta_{\scr E}$ in $K(\P(\scr E\oplus\scr O_\frak X))$ can be written as
\[
\sum_{i=0}^{r}c_{r-i}(\scr E)\boxtimes \scr O(-i),\quad\text{where}\quad c_{r-i}(\scr E) = (-1)^{r-i}{\textstyle(\det(\scr E)^{\vee}\tens \bigwedge^i\scr E)},
\]
and it is trivialized in $K(\P(\scr E))$ via the Koszul complex of the tautological epimorphism $\scr E_{\P(\scr E)}(-1)\onto \scr O_{\P(\scr E)}$.

\begin{definition}
	A $K_{\frak X}$-module is called \emph{Bott periodic} if it is $\beta_{\scr E}$-periodic for every locally free module of finite rank $\scr E$ over $\frak X$.
\end{definition}

In the diagram~\eqref{eqn:PBF}, we can replace $K_\frak X$ by $K^B_\frak X$ \cite[Theorem 3.12 (3)]{KrishnaRavi}, and also by $L_{\A^1}K_\frak X$ or $L_{\A^1}K^B_\frak X$, as the projective bundle formula obviously persists after applying the naive $\A^1$-localization. As a result, all these $K_\frak X$-modules are Bott periodic.

We denote by $\KH_\frak X$ the reflection of $K_\frak X$ in the $\infty$-category of homotopy invariant, Nisnevich excisive, and Bott periodic $K_\frak X$-modules. Since motivic localization and periodization are both compatible with the symmetric monoidal structure, $\KH_\frak X$ is an $E_\infty$-algebra under $K_\frak X$.

\begin{definition}
	The \emph{homotopy $K$-theory} of $\frak X\in\Tame_B$ is the $E_\infty$-ring spectrum $\KH(\frak X)=\KH_\frak X(\frak X)$.
\end{definition}

Property (1) of Theorem~\ref{thm:intro} is clear: if $\frak X$ is regular, $K$-theory is already a homotopy invariant Nisnevich sheaf on $\Sm_\frak X$ \cite[Theorems 2.7, 4.1, and 5.7]{ThomasonK}, and both $L_\mot$ and periodization commute with restriction along $\Sm_\frak X\into \Sch_\frak X$ (since right Kan extension preserves the corresponding local objects).

For any N-quasi-projective morphism $f\colon \frak Y\to \frak X$, let $f^*(\KH_\frak X)$ denote the restriction of $\KH_\frak X$ to $\Sch_\frak Y$.
Then $f^*(\KH_\frak X)$ is the reflection of $K_\frak Y$ in the $\infty$-category of $K_\frak Y$-modules that are homotopy invariant, Nisnevich excisive, and periodic with respect to the maps $\beta_{f^*(\scr E)}$, where $\scr E$ is a locally free module over $\frak X$. In particular, there is a canonical morphism of $E_\infty$-algebras $f^*(\KH_\frak X)\to \KH_\frak Y$.

\begin{proposition}\label{prop:globalBC}
	Let $f\colon\frak Y\to\frak X$ be an N-quasi-projective morphism in $\Tame_B$. Then the map $f^*(\KH_\frak X)\to\KH_\frak Y$ is an equivalence. In other words, $\KH_\frak X$ is the restriction of $\KH$ to $\Sch_\frak X$. 
\end{proposition}

Proposition~\ref{prop:globalBC} immediately implies properties (2) and (3) of Theorem~\ref{thm:intro}, and also that $\KH$ is a Nisnevich sheaf.
Before proving it, we relate the periodization process in the definition of $\KH_\frak X$ to the Bass construction, which will also lead to a proof of property (4).

Consider the map $\beta_{\scr O}\colon \P^1/\infty\to K_\frak X$. The image of $\beta_{\scr O}$ in $K_0(\P^1)$ is thus $[\scr O(-1)]-[\scr O]$. As an element of $\tilde K_0(\P^1,\infty)$, $\beta_{\scr O}$ is determined by any choice of trivialization of $\scr O(-1)$ over $\infty$; we choose the point $(1,0)$ in the line $[1:0]$. This trivialization extends to the standard trivialization of $\scr O(-1)$ over $\P^1\minus 0$, which defines a lift of $\beta_{\scr O}$ to $\P^1/(\P^1\minus 0)$. Moreover, the restriction of the latter to $\A^1=\P^1\minus\infty\subset\P^1$ is nullhomotopic via the standard trivialization of $\scr O(-1)$ over $\P^1\minus\infty$. Since these two trivializations coincide over $1\in\G_m$, we obtain the following commutative diagram of pointed presheaves:
\begin{tikzequation}
	\label{eqn:beta}
	\diagram{\A^1/\G_m & \P^1/(\P^1\minus 0) & \P^1/\infty \\ \Sigma(\G_m/1) & & K_\frak X\rlap. \\};
	\arrows (11-) edge[c->] (-12) (12-) edge[<<-] (-13) (11) edge (21) (21-) edge[dashed] node[below]{$\beta$} (-23) (13) edge node[right]{$\beta_{\scr O}$} (23) (12) edge[dashed] (23);
\end{tikzequation}
The homotopy class of the lower map is the Bott element $\beta\in\tilde K_1(\G_m,1)$ of Example~\ref{ex:KB} (this identification depends on a choice of orientation of the loop in $\Sigma(\G_m/1)$: if the left vertical arrow in~\eqref{eqn:beta} is $*\sqcup_{\G_m}\A^1\to *\sqcup_{\G_m}*$, we let the loop go from the first to the second vertex).
Recall from Example~\ref{ex:KB} that $\gamma$ is the composition of the collapse map $(\P^1\minus 0)\coprod_{\G_m}\A^1\to \Sigma(\G_m/1)$ and $\beta$.

\begin{lemma}\label{lem:Bperiodicity}
	Let $E$ be a $K_\frak X$-module.
	\begin{enumerate}
		\item Suppose that $E$ is $\A^1$-invariant. Then $E$ is $\beta$-periodic if and only if it is $\gamma$-periodic.
		\item Suppose that $E$ is a Zariski sheaf. Then $E$ is $\beta_{\scr O}$-periodic if and only if it is $\gamma$-periodic.
	\end{enumerate}
\end{lemma}

\begin{proof}
	Assertion (1) follows from the fact that $(\P^1\minus 0)\coprod_{\G_m}\A^1\to \Sigma(\G_m/1)$ is an $L_{\A^1}$-equivalence. By~\eqref{eqn:beta}, we can identify $\gamma$ with the composition
	\[
	(\P^1\minus 0)\coprod_{\G_m}\A^1 \to \P^1/1 \onto \P^1/(\P^1\minus 0) \to K_{\frak X},
	\]
	where the first map is a Zariski equivalence.
	Let $\phi\colon \P^1/1\to \P^1/\infty$ be the linear automorphism of $\P^1$ that fixes $0$ and exchanges $1$ and $\infty$.
	Then the square
	\begin{tikzmath}
		\diagram{\P^1/1 & \P^1/\infty \\ \P^1/(\P^1\minus 0) & K_{\frak X} \\};
		\arrows (11-) edge node[above]{$\phi$} (-12) (11) edge[->>] (21) (21-) edge (-22) (12) edge node[right]{$\beta_{\scr O}$} (22);
	\end{tikzmath}
	commutes up to homotopy, since both compositions classify the same element in $\tilde K_0(\P^1,1)$. Assertion (2) follows. 
\end{proof}

Recall from Example~\ref{ex:KB} that $K^B_\frak X=P_\gamma K_\frak X$.
It follows from Lemma~\ref{lem:Bperiodicity} that $\KH_\frak X$ is $\gamma$-periodic (as well as $\beta$-periodic).
Hence, we have morphisms of $E_\infty$-algebras
\[
K_\frak X \to K^B_\frak X \to L_{\A^1}K^B_\frak X \to \KH_\frak X.
\]
By \cite[Theorem 3.12]{KrishnaRavi}, $K^B_\frak X$ is Nisnevich excisive and Bott periodic, so $K^B_\frak X$ is in fact the reflection of $K_\frak X$ in the subcategory of Nisnevich excisive Bott periodic $K_\frak X$-modules. 
Similarly, $L_{\A^1}K^B_\frak X$ is the reflection of $K_\frak X$ in the subcategory of $\A^1$-invariant, Nisnevich excisive, and Bott periodic $K_\frak X$-modules. 
If $\frak X\in\Sch_{\mathbf BG}$ where $G$ is an extension of a finite group scheme by a Nisnevich-locally diagonalizable group scheme, every $\A^1$-invariant Nisnevich sheaf on $\Sch_\frak X$ is already homotopy invariant \cite[Remark 3.13]{Hoyois}, and so the map $L_{\A^1}K^B_\frak X \to \KH_\frak X$ is an equivalence. This proves property (4) of Theorem~\ref{thm:intro}.

We observe that the assignment $\scr E\mapsto\beta_{\scr E}$ is a functor from the groupoid of locally free modules of finite rank over $\frak X$ to the overcategory of $K_\frak X$.
This functoriality comes from~\eqref{eqn:Perf} and the fact that the sheaf $\scr O(-i)$ on $\P(\scr E)$, as $\scr E$ varies in this groupoid, is a cartesian section of the fibered category of quasi-coherent sheaves. 
In particular, $\beta_{\scr E}\colon\V^+(\scr E)\to K_\frak X$ coequalizes the action of linear automorphisms of $\scr E$ on $\V^+(\scr E)$.

Write $\V_0(\scr E)$ and $\V_0^+(\scr E)$ for the pointed presheaves $\V(\scr E)/(\V(\scr E)\minus 0)$ and  $\P(\scr E\oplus\scr O_\frak X)/(\P(\scr E\oplus\scr O_\frak X)\minus 0)$ on $\Sch_\frak X$. 
As in~\eqref{eqn:beta}, we have a zig-zag
\[
\V_0(\scr E)\into \V_0^+(\scr E) \twoheadleftarrow \V^+(\scr E),
\]
where the first map is a Zariski equivalence and the second map is an $L_{\A^1}$-equivalence.
The map $\beta_{\scr E}\colon \V^+(\scr E)\to K_\frak X$ extends to $\V_0^+(\scr E)$ because the morphism~\eqref{eqn:Koszul} is an epimorphism away from the zero section, and hence it induces
\[
\beta'_\scr E\colon \V_0(\scr E)\to K_\frak X.
\]
Explicitly, $\beta'_\scr E$ is represented by the Koszul complex of the tautological morphism $\scr E_{\V(\scr E)}\to \scr O_{\V(\scr E)}$ tensored with $\det(\scr E)[r]^\vee$, viewed as an object of $\Perf(\V(\scr E)\on \frak X)$.

Note that the assignment $\scr E\mapsto \V_0(\scr E)$ is right-lax symmetric monoidal, with the monoidal structure maps $\V_0(\scr E)\otimes \V_0(\scr F)\to \V_0(\scr E\oplus\scr F)$ being Zariski equivalences. Using the Koszul complex representative of $\beta_\scr E'$ and the multiplicative properties of Koszul complexes, we can promote the assignment $\scr E\mapsto\beta'_\scr E$ to a right-lax symmetric monoidal functor from the groupoid of locally free modules of finite rank over $\frak X$ (under direct sum) to the $\infty$-category of presheaves of spectra on $\Sch_\frak X$ over $K_\frak X$.
In particular, if $\scr E$ and $\scr F$ are locally free modules of finite rank over $\frak X$, we have a commutative square
\begin{tikzmath}
	\def\colsep{5em}
	\diagram{\V_0(\scr E)\otimes \V_0(\scr F) & K_\frak X\otimes K_\frak X \\ \V_0(\scr E\oplus\scr F) & K_{\frak X}\rlap, \\};
	\arrows (11-) edge node[above]{$\beta_\scr E'\otimes\beta_\scr F'$} (-12) (11) edge (21) (21-) edge node[above]{$\beta_{\scr E\oplus\scr F}'$} (-22) (12) edge (22);
\end{tikzmath}
where the right vertical map is multiplication.

\begin{proof}[Proof of Proposition~\ref{prop:globalBC}]
	We must show that the map $f^*(\KH_\frak X)\to \Hom(\V^+(\scr E),f^*(\KH_\frak X))$ induced by $\beta_\scr E$ is an equivalence for every locally free module $\scr E$ over $\frak Y$.
	Since $f^*(\KH_\frak X)$ is a homotopy invariant Nisnevich sheaf, we can assume that $\frak Y\to\frak X$ is quasi-affine \cite[Proposition 4.6]{Hoyois}.
	By Lemma~\ref{lem:complete}, we can then write $\scr E$ as a quotient of $f^*(\scr G)$ for some locally free module of finite rank $\scr G$ over $\frak X$. Replacing $\frak Y$ by an appropriate vector bundle torsor, we can assume that $f^*(\scr G)\simeq \scr E\oplus\scr F$ for some $\scr F$.
	Hence, $\beta'_{\scr E\oplus\scr F}\simeq\beta'_{\scr E}\beta'_{\scr F}$ acts invertibly on $f^*(\KH_\frak X)$. 
	In the sequence
	\[
		f^*(\KH_\frak X)\xrightarrow{\beta'_\scr E} \Hom(\V_0(\scr E),f^*(\KH_\frak X)) \xrightarrow{\beta'_\scr F} \Hom(\V_0(\scr E\oplus\scr F),f^*(\KH_\frak X)) \xrightarrow{\beta'_\scr E} \Hom(\V_0(\scr E\oplus\scr F\oplus\scr E),f^*(\KH_\frak X)),
	\]
	the composites $\beta'_\scr F\beta'_\scr E$ and $\beta'_\scr E\beta'_\scr F$ are thus both equivalences. It then follows from the 2-out-of-6 property that all three maps are equivalences.
\end{proof}

This concludes the verification of properties (1)--(4) of Theorem~\ref{thm:intro}.
Finally, we would like to obtain a more concrete description of $\KH$ using Theorem~\ref{thm:P=Q}.
In the following lemma, $\Sp(\PSh_{\A^1,\Zar}(\Sch_\frak X))$ denotes the $\infty$-category of $\A^1$-invariant Zariski sheaves of spectra on $\Sch_\frak X$.
 
\begin{lemma}\label{lem:Bottsymmetric}
	Let $\scr E$ be a locally free module of finite rank over $\frak X$.
	Then $L_{\A^1,\Zar}\beta_{\scr E}'\colon L_{\A^1,\Zar}\V_0(\scr E)\to L_{\A^1,\Zar} K_\frak X$, viewed as an object of $\Sp(\PSh_{{\A^1,\Zar}}(\Sch_\frak X))_{/L_{\A^1,\Zar} K_\frak X}$, is $3$-symmetric.
\end{lemma}

\begin{proof}
	Since $\scr E\mapsto L_{\A^1,\Zar} \beta'_\scr E$ is symmetric monoidal, it suffices to show that $L_{\A^1}\sigma_3\colon L_{\A^1}\V_0(\scr E^3)\to L_{\A^1}\V_0(\scr E^3)$ is homotopic to the identity over $L_{\A^1} K_\frak X$.
	 The identity and $\sigma_3$ are both induced by matrices in $\mathrm{SL}_3(\Z)$ acting on $\scr E^3$, and any two such matrices are $\A^1$-homotopic. Thus, it will suffice to prove the following statement: for any locally free module of finite rank $\scr E$ over $\frak X$ and any automorphism $\phi$ of $p^*(\scr E)$, where $p\colon\A^1\times\frak X\to\frak X$ is the projection, the automorphisms of $\V_0(\scr E)$ induced by $\phi_0$ and $\phi_1$ are $\A^1$-homotopic over $L_{\A^1}K_\frak X$.
	 Since $\beta'_{\scr E}$ is functorial in $\scr E$, the automorphism $\phi$ induces a commutative triangle
	  \begin{tikzmath}
		\diagram[column sep={4em,between origins}]{\V_0(p^*(\scr E)) & & \V_0(p^*(\scr E)) \\ & L_{\A^1}K_{\A^1\times\frak X} & \\};
		\arrows (11-) edge node[above]{$\phi$} (-13) (11) edge node[left]{$\beta'_{p^*(\scr E)}$} (22) (13) edge node[right]{$\beta'_{p^*(\scr E)}$} (22);
	  \end{tikzmath}
	  of presheaves of spectra on $\Sch_{\A^1\times\frak X}$.
	  By adjunction, this is equivalent to a triangle
	  \begin{tikzmath}
		\diagram[column sep={4em,between origins}]{\A^1_+\otimes \V_0(\scr E) & & \V_0(\scr E) \\ & L_{\A^1}K_{\frak X}\rlap, & \\};
		\arrows (11-) edge (-13) (11) edge node[left]{$\beta'_{\scr E}$} (22) (13) edge node[right]{$\beta'_{\scr E}$} (22);
	  \end{tikzmath}
	  which is an $\A^1$-homotopy between $\phi_0$ and $\phi_1$ over $L_{\A^1}K_\frak X$, as desired.
\end{proof}

\begin{proposition}\label{prop:KHconcrete}
	Let $\frak X\in\Tame_B$ and let $E$ be a $K_\frak X$-module. Then the canonical map $E\to Q_{\{\beta_\scr E\}} L_\mot E$ is the universal map to a homotopy invariant, Nisnevich excisive, and Bott periodic $K_\frak X$-module. In particular,
	\[
	\KH_\frak X \simeq Q_{\{\beta_\scr E\}} L_\mot K_\frak X.
	\]
\end{proposition}

\begin{proof}
	Combining Lemma~\ref{lem:Bottsymmetric} and Theorem~\ref{thm:P=Q}, we deduce that
	\[
	P_{\{\beta'_\scr E\}}L_\mot E\simeq Q_{\{\beta'_\scr E\}}L_\mot E.
	\]
	As $L_\mot E$ is in particular a Zariski sheaf, we can replace $\beta'_\scr E$ with $\beta_\scr E$ without changing either side.
	Hence, we have
	\[
	P_{\{\beta_\scr E\}} L_\mot E \simeq Q_{\{\beta_\scr E\}} L_\mot E.
	\]
	We conclude by noting that $Q_{\{\beta_\scr E\}}$ preserves homotopy invariant Nisnevich sheaves.
\end{proof}

In other words, $\KH_\frak X$ is the Bott spectrification of the motivic localization of $K_\frak X$.

\section{The equivariant motivic $K$-theory spectrum}
\label{sec:cdh}

In this final section, we prove that $\KH$ is a cdh sheaf on $\Tame_B$. 
By definition of the cdh topology, this is the case if and only if the restriction of $\KH$ to $\Sch_\frak X$ is a cdh sheaf for every $\frak X\in\Tame_B$. Moreover, as we already know that $\KH$ is a Nisnevich sheaf, we can assume without loss of generality that $\frak X=[X/G]$ with $X$ a small $G$-scheme. By definition of smallness, we may as well assume that $B$ has the $G$-resolution property and that $\frak X=\mathbf BG$. Thus, we are now in the setting of \cite[\S6]{Hoyois}.

Let $\H_\pt(\frak X)$ be the pointed motivic homotopy category over $\frak X\in\Sch_\BG$, \ie, the $\infty$-category of pointed presheaves on $\Sm_\frak X$ that are homotopy invariant and Nisnevich excisive. The stable motivic homotopy category over $\frak X$ is by definition
\[
\SH(\frak X) = \H_\pt(\frak X)[\Sph_\BG^{-1}],
\]
where $\Sph_\BG$ is the collection of one-point compactifications $\V^+(\scr E)$ of vector bundles over $\BG$ (pulled back to $\frak X$); this forces the invertibility of the one-point compactifications of \emph{all} vector bundles over $\frak X$ \cite[Corollary 6.7]{Hoyois}.
Let $\Sp(\H(\frak X))$ be the $\infty$-category of homotopy invariant Nisnevich sheaves of spectra on $\Sm_\frak X$, or equivalently the stabilization of $\H(\frak X)$. As a symmetric monoidal $\infty$-category, it is $\H_\pt(\frak X)[(S^1)^{-1}]$. Since $S^1$ is invertible in $\SH(\frak X)$, we have
\[
\SH(\frak X) \simeq \Sp(\H(\frak X))[\Sph_\BG^{-1}].
\]

We also consider ``big'' variants of these $\infty$-categories: $\Sp(\underline\H)$ is the $\infty$-category of homotopy invariant Nisnevich sheaves of spectra on $\Sch_\BG$, and $\underline\SH=\Sp(\underline\H)[\Sph_\BG^{-1}]$.
These $\infty$-categories have the following interpretation. Any presheaf on $\Sch_\BG$ can be restricted to $\Sm_\frak X$ for every $\frak X\in\Sch_\BG$; this gives rise to a section of the cocartesian fibration classified by $\frak X\mapsto \PSh(\Sm_\frak X)$, which sends smooth morphisms to cocartesian edges. It is clear that this construction is an equivalence of $\infty$-categories between presheaves on $\Sch_\BG$ and such sections. From this we deduce that $\Sp(\underline\H)$ and $\underline\SH$ can be identified with $\infty$-categories of sections of $\Sp(\H(\ph))$ and $\SH(\ph)$ over $\Sch_\BG^\op$ that are cocartesian over smooth morphisms.

In \S\ref{sec:KH}, we constructed the $E_\infty$-algebra $\KH_\BG$ in $\Sp(\underline\H)$ as a Bott periodic $K_\BG$-module.
By Proposition~\ref{prop:Pdeloop}, there is a unique Bott periodic $E_\infty$-algebra $\underline\KGL$ in $\underline\SH$ such that $\Omega^\infty\underline\KGL\simeq \KH_\BG$, namely
\[
\underline\KGL = P_{\{\beta_\scr E\}}\Sigma^\infty\KH_\BG.
\]
By Proposition~\ref{prop:KHconcrete}, we can write $\underline\KGL$ more explicitly as an $\Sph_\BG$-spectrum in $\Sp(\underline\H)$: it is the image, under the localization functor
\[
QL_\mot\colon \Stab^\lax_{\Sph_\BG}\Sp(\PSh(\Sch_\BG)) \to \Stab_{\Sph_\BG}\Sp(\underline\H)\simeq\underline\SH,
\]
of the ``constant'' $\Sph_\BG$-spectrum $c_{\{\beta_\scr E\}} K_\BG$.

\begin{definition}
	For $\frak X\in\Sch_\BG$, we denote by $\KGL_\frak X\in\CAlg(\SH(\frak X))$ the restriction of $\underline\KGL$ to $\Sm_\frak X$.
\end{definition}

By Proposition~\ref{prop:globalBC}, the motivic spectrum $\KGL_\frak X$ represents homotopy $K$-theory: for $\frak Y$ a smooth N-quasi-projective $\frak X$-stack, there is a natural equivalence
\[
\KH(\frak Y) \simeq \Mor(\Sigma^\infty\frak Y_+, \KGL_\frak X),
\]
where $\Mor$ denotes a mapping spectrum in the stable $\infty$-category $\SH(\frak X)$.

We now prove that $\frak X\mapsto\KGL_\frak X$ is a cocartesian section of $\SH(\ph)$ over $\Sch_\BG^\op$, \ie, that for every $f\colon\frak Y\to\frak X$ in $\Sch_\BG$, the restriction map
\[
f^*(\KGL_\frak X)\to\KGL_\frak Y
\]
in $\SH(\frak Y)$ is an equivalence. By \cite[Corollary 6.25]{Hoyois}, this implies that $\KH$ is a cdh sheaf on $\Sch_\BG$ and concludes the proof of Theorem~\ref{thm:intro}.
Since $\underline\KGL=QL_\mot c_{\{\beta_\scr E\}}K_\BG$, the above restriction map is
\[
f^*(QL_\mot c_{\{\beta_\scr E\}}(K|\Sm_\frak X))\to QL_\mot c_{\{\beta_\scr E\}}(K|\Sm_\frak Y)
\]
The localization functor $QL_\mot$ is compatible with the base change functor $f^*$, as $f_*$ preserves local objects, so it will suffice to show that the restriction map
\begin{equation}\label{eqn:resK}
	f^*(K|\Sm_\frak X) \to K|\Sm_\frak Y
\end{equation}
is a motivic equivalence in $\Sp(\PSh(\Sm_\frak Y))$.

Sending vector bundles over $\frak X$ to their classes in $K$-theory induces a map of grouplike $E_\infty$-spaces
\begin{equation}\label{eqn:group_completion}
	\Vect(\frak X)^+\to \Omega^\infty K(\frak X),
\end{equation}
where $\Vect(\frak X)$ is the $E_\infty$-space of vector bundles over $\frak X$ and $(\ph)^+$ denotes group completion.
If $\frak X=[X/G]$ with $X$ a small affine $G$-scheme, it follows from \cite[Lemma 2.17]{Hoyois} that every short exact sequence of vector bundles over $\frak X$ splits. In that case, the map~\eqref{eqn:group_completion} is an equivalence.
By \cite[Proposition 3.16 (2)]{Hoyois}, it follows that the map
\[
\Vect^+\to \Omega^\infty K|\Sm_\frak X
\]
is a motivic equivalence in $\PSh(\Sm_\frak X)$. Note also that the inclusion
\[
\coprod_{n\geq 0}B_\fppf \GL_n\into\Vect
\]
exhibits $\Vect$ as the Zariski sheafification of the subgroupoid of vector bundles of constant rank. By Lemma~\ref{lem:plus} below, it remains a Zariski equivalence after group completion. We therefore obtain a motivic equivalence
\begin{equation}\label{eqn:BGL_nK}
\Biggl(\coprod_{n\geq 0}B_\fppf \GL_n\Biggr)^+\to \Omega^\infty K|\Sm_\frak X.
\end{equation}

\begin{lemma}\label{lem:plus}
	Let $F\colon \scr C\to\scr D$ be a colimit-preserving functor between presentable $\infty$-categories. Suppose that finite products distribute over colimits in $\scr C$ and $\scr D$ and that $F$ preserves finite products. Then, for every $E_\infty$-monoid $M$ in $\scr C$, the canonical map $F(M)^+\to F(M^+)$ is an equivalence.
\end{lemma}

\begin{proof}
	The assumption on $\scr C$ implies that the $\infty$-category $\CMon(\scr C)$ of $E_\infty$-monoids in $\scr C$ is presentable \cite[Corollary 3.2.3.5]{HA} and hence that group completion exists.
	Since both $F$ and its right adjoint preserve finite products, they lift to a pair of adjoint functors between $\CMon(\scr C)$ and $\CMon(\scr D)$, as well as between the subcategories of grouplike objects. This immediately implies that $F$ commutes with group completion.
\end{proof}

For any $f\colon\frak Y\to\frak X$ in $\Sch_\BG$, the pullback functor $f^*\colon\PSh(\Sm_\frak X)\to\PSh(\Sm_\frak Y)$ preserves finite products and hence commutes with group completion of $E_\infty$-monoids, by Lemma~\ref{lem:plus}.
Similarly, since $L_\mot\colon \PSh(\Sm_\frak X)\to \H(\frak X)$ preserves finite products \cite[Proposition 3.15]{Hoyois}, it commutes with group completion of $E_\infty$-monoids.
Hence, by~\eqref{eqn:BGL_nK} and Corollary~\ref{cor:BC} (with $\Gamma=\GL_n$), we deduce that the restriction map
\[
f^*(\Omega^\infty K|\Sm_\frak X) \to \Omega^\infty K|\Sm_\frak Y
\]
is a motivic equivalence in the $\infty$-category of grouplike $E_\infty$-monoids in $\PSh(\Sm_\frak Y)$.
Equivalently, \eqref{eqn:resK} is a motivic equivalence in $\Sp_{\geq 0}(\PSh(\Sm_\frak Y))$, whence in $\Sp(\PSh(\Sm_\frak Y))$, as was to be shown.
 
 \begin{remark}
 	If $f\colon\frak Y\to \frak X$ is a morphism of \emph{schemes}, it is easy to show that the map~\eqref{eqn:resK} is a Zariski equivalence, because $B_\fppf\GL_n=B_\Zar\GL_n$ and $\GL_n$ is smooth. The proof of cdh descent in this case does not need the geometric model for the classifying space of $\GL_n$. 
 \end{remark}
 
 \begin{proof}[Comments on Theorem~\ref{thm:intro-2}]
 	We discuss the minor modifications needed for the proof of Theorem~\ref{thm:intro-2}.
	If $X$ is a locally affine qcs $G$-schemes such that $\abs G$ is invertible on $X$, then $[X/G]$ is a qcs tame Deligne–Mumford stack with coarse moduli scheme. By \cite[Corollary 3.8]{KO} and a noetherian approximation argument, nonconnective $K$-theory is a Nisnevich sheaf on such stacks, whence also $\KH$ (defined as the naive $\A^1$-localization of $K^B$). The projective bundle formula holds for general stacks \cite[Theorem 3.6]{KrishnaRavi}. Hence, the restriction of $\KH$ to the category of smooth quasi-projective $G$-schemes over $X$ is a homotopy invariant Nisnevich sheaf as well as a Bott periodic $E_\infty$-algebra. By Proposition~\ref{prop:Pdeloop}, it deloops uniquely to a Bott periodic $E_\infty$-algebra $\KGL_{[X/G]}\in\SH([X/G])$. Since $[X/G]$ is Nisnevich-locally of the form $[U/G]$ with $U$ affine, the proof of Theorem~\ref{thm:MV} and the above arguments go through (with some simplifications) and show that, for every $G$-equivariant morphism $f\colon Y\to X$ with $Y$ a locally affine qcs $G$-scheme, $f^*(\KGL_{[X/G]})\simeq \KGL_{[Y/G]}$. By \cite[Remark 6.26]{Hoyois}, we conclude that $\KH$ satisfies cdh descent on the category of locally affine qcs $G$-schemes.
 \end{proof}
 
 \begin{proof}[Comments on Theorem~\ref{thm:intro-3}]
 	Because of the reductions done at the beginning of this section, we have only proved Theorem~\ref{thm:intro-3} with $\Tame_B$ replaced by the subcategory of stacks $\frak X$ admitting an N-quasi-projective map $\frak X\to\mathbf B_UG$ for some $B$-scheme $U$ such that $\mathbf B_UG$ has the resolution property. In fact, $\SH(\frak X)$ is only defined for such $\frak X$ in \cite[\S6]{Hoyois}. As indicated in \emph{loc.\ cit.}, however, $\SH(\ph)$ extends uniquely, by right Kan extension, to a Nisnevich sheaf on $\Tame_B$. Hence, the section $\frak X\mapsto\KGL_\frak X$ constructed above also extends uniquely to a section of $\CAlg(\SH(\ph))$ on all of $\Tame_B^\op$ that is cocartesian over N-quasi-projective morphisms, and Theorem~\ref{thm:intro-3} holds in the stated generality.
 \end{proof}

\begin{remark}
	Suppose that $\frak X\in\Tame_B$ is regular noetherian. Then the Borel–Moore homology theory on $\Sch_\frak X$ represented by $\KGL_\frak X$ is the $K$-theory of coherent sheaves. More precisely, for every quasi-projective morphism $f\colon\frak Z\to\frak X$, there is an equivalence of spectra
	\[
	\Mor(\1_\frak Z, f^!\KGL_\frak X) \simeq K(\mathrm{Coh}(\frak Z)),
	\]
	where the left-hand side is a mapping spectrum in $\SH(\frak Z)$. To prove this, write $f=p\circ i$ where $i\colon \frak Z\into \frak Y$ is a closed immersion and $p\colon\frak Y\to\frak X$ is smooth quasi-projective. 
	By \cite[Theorem 6.18 (2)]{Hoyois} and Bott periodicity, $p^!\KGL_\frak X \simeq \Sigma^{\Omega_p}\KGL_\frak Y\simeq \KGL_\frak Y$.
	Let $j\colon\frak U\into \frak Y$ be the open immersion complementary to $i$. By \cite[Theorem 6.18 (4)]{Hoyois}, we have a fiber sequence
	\[
	i_*i^!\KGL_\frak Y \to \KGL_\frak Y\to j_*j^*\KGL_\frak Y
	\]
	in $\SH(\frak Y)$, whence a fiber sequence of spectra
	\[
	\Mor(\1_\frak Z,f^!\KGL_\frak X) \to \Mor(\1_{\frak Y},\KGL_\frak Y) \to \Mor(\1_\frak U,\KGL_\frak U).
	\]
	Since $\frak Y$ and $\frak U$ are regular, the second map is identified with the restriction $K(\frak Y)\to K(\frak U)$, whose fiber is $K(\mathrm{Coh}(\frak Z))$.
\end{remark}

\providecommand{\bysame}{\leavevmode\hbox to3em{\hrulefill}\thinspace}

\end{document}